\documentclass[10pt]{article}
\usepackage{amsfonts}
\usepackage{hyperref}
\usepackage{mathrsfs,enumerate,amssymb,amsmath,color,lineno}
\usepackage{indentfirst}
\usepackage{amsmath}
\usepackage{amssymb}
\usepackage[amsmath,thmmarks]{ntheorem}
\usepackage{graphicx}
\usepackage{tikz}

\newtheorem{definition}{\bf Definition}[section]
\newtheorem{lemma}{\bf Lemma}[section]
\newtheorem{theorem}{\bf Theorem}[section]
\newtheorem{remark}{\bf Remark}[section]
\newtheorem{corollary}{\bf Corollary}[section]
\newtheorem{example}{\bf Example}[section]

\newtheorem{proposition}{\bf Proposition}[section]

\begin{document}
\setcounter{page}{1}
	
\title{{\textbf{Additive generator pairs of overlap functions}}\thanks {Supported by the National Natural Science Foundation of China (No.12471440)}}
\author{Li-zhi Liang\footnote{\emph{E-mail address}: 1902428258@qq.com}, Xue-ping Wang\footnote{Corresponding author. xpwang1@hotmail.com; fax: +86-28-84761502},\\
\emph{School of Mathematical Sciences, Sichuan Normal University,}\\
\emph{Chengdu 610066, Sichuan, People's Republic of China}}
	
\newcommand{\pp}[2]{\frac{\partial #1}{\partial #2}}
\date{}
\maketitle
\begin{quote}
{\bf Abstract} Let $\theta:[0,1]\rightarrow[-\infty,+\infty]$ be a function with both $\theta(x^{-})$ and $\theta(x^{+})$ existing for every $x\in [0,1]$ and $\vartheta:[-\infty,+\infty]\rightarrow[-\infty,+\infty]$ be a function. In this article we completely characterize the pair $(\theta,\vartheta)$ for the bivariate function $O_{\theta,\vartheta}:[0,1]^{2}\rightarrow[0,1]$ given by
	$$O_{\theta,\vartheta}(x,y)=\vartheta(\theta(x)+\theta(y))$$
	being an overlap function. In particular, we give analytical expressions of some transformations for the pair $(\theta,\vartheta)$.
		
\end{quote}
	
\section{Introduction}
Although triangular norms and triangular conorms are widely applied to these fields such as image processing, fuzzy logic, classification and decision making in the real world \cite{JF1994}, their associativity is strongly required. To avoid this problem, arising from some practical applications in image processing, classification and decision making Bustince et al. introduced overlap and grouping functions \cite{HB2010,HB2012}. After that, overlap and grouping functions undergo a rapid development in two aspects included applications and theory. In applications, overlap and grouping functions are widely used to solve difficulties appeared in image processing \cite{HB2010}, data stream clustering \cite{AU2025}, convolutional neural networks \cite{IR2023}, classification \cite{TD2021,LI2025} and decision making problems \cite{RS2021}. Accompanying with the development of applications, scholars all over the world also obtain remarkable progress in theory of overlap and grouping functions, such as their migrativity, cross-migrativitity, bimigrativity, homogeneity, idempotency, the Lipschitz condition and the existence of generators \cite{BB2013,BB2014,GP2015,GP2014,FA2025,DG2016,WA2025}.

It is well-known that additive and multiplicative generators play a significant role for building t-norms or t-conorms \cite{EP2000,OY2007,PV2005}, so that it seems natural for us to consider the use of additive or multiplicative generators for the overlap and grouping functions. Moreover, the use of additive or multiplicative generators may simplify the choice of an appropriate overlap or grouping function for a given problem since we only need to consider one-variable functions instead of a bivariate one, reducing the computational complexity in this way \cite{AM2009}. Following this idea, Dimuro et al. \cite{GD2014} introduced the notion of an additive generator pair for a grouping function. Furthermore, they gave the concept of an additive generator pair for an overlap function and explored the overlap functions generated by additive generator pairs. Specific speaking, they gave the definition of an additive generator pair for an overlap function as follows.
\begin{definition}[\cite{GP2016}]\label{def001}
\emph{Let $\theta:[0,1]\rightarrow[0,+\infty]$ and $\vartheta:[0,+\infty]\rightarrow[0,1]$ be continuous and non-increasing functions such that
    \item[(1)] $\theta(x)=+\infty$ if and only if $x=0$;
	\item[(2)] $\theta(x)=0$ if and only if $x=1$;
	\item[(3)] $\vartheta(x)=1$ if and only if $x=0$;
	\item[(4)] $\vartheta(x)=0$ if and only if $x=+\infty$.\\
Then the function $O_{\theta,\vartheta}:[0,1]^{2}\rightarrow[0,1]$, defined by $O_{\theta,\vartheta}(x,y)=\vartheta(\theta(x)+\theta(y))$, is an overlap function. $(\theta,\vartheta)$ is called an additive generator pair of the overlap function $O_{\theta,\vartheta}$, and $O_{\theta,\vartheta}$ is said to be additively generated by the pair $(\theta,\vartheta)$.}	
\end{definition}

Based on Definition \ref{def001}, we recently revealed the relationship between functions $\theta$ and $\vartheta$ in an overlap function additively generated by an additive generator pair ($\theta$,$\vartheta$), which is used to characterize the conditions for an overlap function additively generated by the pair ($\theta$,$\vartheta$) being a triangular norm by terms of functions $\theta$ and $\vartheta$ \cite{LL2025}.

From Definition \ref{def001}, a pair $(\theta,\vartheta)$ which satisfies the conditions (1), (2), (3), (4) and the function $O_{\theta,\vartheta}:[0,1]^{2}\rightarrow[0,1]$ defined by $O_{\theta,\vartheta}(x,y)=\vartheta(\theta(x)+\theta(y))$ is an overlap function is called an additive generator pair of the overlap function $O_{\theta,\vartheta}$. Thereby, an interesting question arises: what does a pair $(\theta,\vartheta)$ need to satisfy whenever the function $O_{\theta,\vartheta}:[0,1]^{2}\rightarrow[0,1]$ defined by $O_{\theta,\vartheta}(x,y)=\vartheta(\theta(x)+\theta(y))$ is an overlap function?

In this article, we directly give the definition of an additive generator pair of an overlap function as follows.
\begin{definition}
\emph{Let $\theta:[0,1]\rightarrow[0,+\infty]$ and $\vartheta:[0,+\infty]\rightarrow[-\infty,+\infty]$ be two functions. If the bivariate function $O_{\theta,\vartheta}:[0,1]^{2}\rightarrow[0,1]$ given by
$$O_{\theta,\vartheta}(x,y)=\vartheta(\theta(x)+\theta(y))$$
is an overlap function, then $(\theta,\vartheta)$ is called an additive generator pair of the overlap function $O_{\theta,\vartheta}$ and $O_{\theta,\vartheta}$ is said to be additively generated by the pair $(\theta,\vartheta)$.}	
\end{definition}
And we positively answer the above problem.

The rest of this article is organized as follows: In Section \ref{section3}, we characterize an additive generator pair $(\theta,\vartheta)$ for the bivariate function $O_{\theta,\vartheta}:[0,1]^{2}\rightarrow[0,1]$ given by
	$$O_{\theta,\vartheta}(x,y)=\vartheta(\theta(x)+\theta(y))$$
	being an overlap function when $\theta:[0,1]\rightarrow[0,+\infty]$ is a function with both $\theta(x^{-})$ and $\theta(x^{+})$ existing for every $x\in [0,1]$ and $\vartheta:[0,+\infty]\rightarrow[-\infty,+\infty]$ is a function. In Section \ref{section4}, we give analytical expressions of some transformations for the pair $(\theta,\vartheta)$ and further describe an additive generator pair $(\theta,\vartheta)$ for the bivariate function $O_{\theta,\vartheta}$ being an overlap function when $\theta:[0,1]\rightarrow[-\infty,+\infty]$ is a function with both $\theta(x^{-})$ and $\theta(x^{+})$ existing for every $x\in [0,1]$ and $\vartheta:[-\infty,+\infty]\rightarrow[-\infty,+\infty]$ is a function. A conclusion is drawn in Section \ref{section5}.


\section{Additive generator pairs of overlap functions}\label{section3}
Let $\theta:[0,1]\rightarrow[0,+\infty]$ be a function with both $\theta(x^{-})$ and $\theta(x^{+})$ existing for every $x\in [0,1]$ and $\vartheta:[0,+\infty]\rightarrow[-\infty,+\infty]$ be a function. In this section we characterize the pair $(\theta,\vartheta)$ for the bivariate function $O_{\theta,\vartheta}:[0,1]^{2}\rightarrow[0,1]$ given by
	$$O_{\theta,\vartheta}(x,y)=\vartheta(\theta(x)+\theta(y))$$
	being an overlap function.

\begin{definition}[\cite{HB2010}]\label{definition2.1}
\emph{A bivariate function $O:[0,1]^{2}\rightarrow[0,1]$ is said to be an overlap function if it satisfies the following conditions:
    \begin{itemize}
	    \item[(O1)] $O$ is commutative,
	    \item[(O2)] $O(x,y)=0$ if and only if $xy=0$,
	    \item[(O3)] $O(x,y)=1$ if and only if $xy=1$,
	    \item[(O4)] $O$ is non-decreasing and
	    \item[(O5)] $O$ is continuous.
    \end{itemize}}
\end{definition}


The following propositions show the conditions what $\theta$ needs to satisfy whenever $O_{\theta,\vartheta}(x,y)=\vartheta(\theta(x)+\theta(y))$ is an overlap function.

\begin{proposition}\label{proposition3.1}
	Let $\theta:[0,1]\rightarrow[0,+\infty]$ be a function with both $\theta(x^{-})$ and $\theta(x^{+})$ existing for every $x\in [0,1]$ and $\vartheta:[0,+\infty]\rightarrow[-\infty,+\infty]$ be a function. If the bivariate function $O_{\theta,\vartheta}:[0,1]^{2}\rightarrow[0,1]$ given by
	$$O_{\theta,\vartheta}(x,y)=\vartheta(\theta(x)+\theta(y))$$
	is an overlap function, then the following two statements hold:
	 \begin{itemize}
	 	\item[(1)] $\theta$ is continuous in $0$;
	 	\item[(2)] $\theta(x)=+\infty$ if and only if $x=0$.
	 \end{itemize}	
\end{proposition}
\begin{proof}
	(1) Supposing that $\theta$ is not continuous in $0$, we have $\theta(0)\ne\theta(0^{+})$ since $\theta:[0,1]\rightarrow[0,+\infty]$ is a function with both $\theta(x^{-})$ and $\theta(x^{+})$ existing for every $x\in [0,1]$. Thus there exists a $u\in(0,1]$ such that $\theta$ is continuous in $(0,u]$. Because $O_{\theta,\vartheta} (x,y)=\vartheta(\theta(x)+\theta(y))$ is an overlap function, it follows from Definition \ref{definition2.1} (O2) and (O5) that $\vartheta(\theta(0^{+})+\theta(u))=\vartheta(\theta(0)+\theta(u))=O_{\theta,\vartheta}(0,u)=0$. Moreover, there is a $w\in(0,u]$ such that $2\theta(w)=\theta(0^{+})+\theta(u)$ since $\theta$ is continuous in $(0,u]$. Thus $O_{\theta,\vartheta}(w,w)=\vartheta(\theta(w)+\theta(w))=\vartheta(\theta(0^{+})+\theta(u))=\vartheta(\theta(0)+\theta(u))=0$ with $w\ne0$, contrary to Definition \ref{definition2.1} (O2). Therefore, $\theta$ is continuous in $0$.
	
	(2)  Assume $\theta(0)\ne+\infty$. Then from (1) there exists a $u\in (0,1]$ such that $\theta$ is continuous in $[0,u]$ and $\theta(u)\ne+\infty$. So there is a $w\in(0,u]$ satisfying $2\theta(w)=\theta(0)+\theta(u)$. Thus $O_{\theta,\vartheta}(w,w)=\vartheta(2\theta(w))=\vartheta(\theta(0)+\theta(u))=0$ with $w\ne0$, contrary to Definition \ref{definition2.1} (O2). Therefore, $\theta(0)=+\infty$.
	
	Conversely, it follows from the above proof that $\theta(0)=+\infty$. If there exists an $x_{0}\in(0,1]$ such that $\theta(x_{0})=+\infty$, then $O_{\theta,\vartheta}(x_0,x_0)=\vartheta(2\theta(x_{0}))=\vartheta(+\infty)=\vartheta(2\theta(0))=0$ with $x_{0}\ne0$, contrary to Definition \ref{definition2.1} (O2).
\end{proof}

From the proof of Proposition \ref{proposition3.1}, the following remark is true.
\begin{remark}
\emph{Under the conditions of Proposition \ref{proposition3.1}, $\vartheta(+\infty)=0$.}
\end{remark}

\begin{proposition}\label{proposition03.2}
	Let $\theta:[0,1]\rightarrow[0,+\infty]$ be a function with both $\theta(x^{-})$ and $\theta(x^{+})$ existing for every $x\in [0,1]$. Then there is a function $\vartheta:[0,+\infty]\rightarrow[-\infty,+\infty]$ such that the bivariate function $O_{\theta,\vartheta}:[0,1]^{2}\rightarrow[0,1]$ given by
	$$O_{\theta,\vartheta}(x,y)=\vartheta(\theta(x)+\theta(y))$$
	is an overlap function if and only if the following four conditions hold:
\begin{itemize}
  \item[(i)] $\theta$ is continuous;
  \item[(ii)] $\theta$ is non-increasing;
  \item[(iii)] $\theta(x)=+\infty$ if and only if $x=0$;
  \item[(iv)] $\theta(x)>\theta(1)$ for each $x\in[0,1)$.
\end{itemize}
	\end{proposition}
\begin{proof}
	(i) If $\theta$ is not continuous, then there is a $u\in [0,1]$ such that $\theta(u^{-})\ne\theta(u^{+})$ and from Proposition \ref{proposition3.1} (1) we have that $\theta$ is continuous in $[0,u)$.	
We just prove that $\theta(u^{-})>\theta(u^{+})$ leads to a contradiction, the case $\theta(u^{-})<\theta(u^{+})$ being analogous. Indeed, because $O_{\theta,\vartheta}(x,y)=\vartheta(\theta(x)+\theta(y))$ is continuous we have
\begin{equation}\label{eq1}
	\vartheta(\theta(x)+\theta(u^{-}))=\vartheta(\theta(x)+\theta(u^{+}))
\end{equation} for each $x\in[0,u)$. Meanwhile, from Proposition \ref{proposition3.1} (2) there exists an $x_{0}\in(0,u]$ such that $\theta$ is non-increasing in $[0,x_{0}]$, and from $\theta(u^{-})>\theta(u^{+})$ we have $\theta(x_{0})+\theta(u^{-})>\theta(x_{0})+\theta(u^{+})$. Thus there is an $x_{1}\in[0,x_{0})$ satisfying
\begin{equation}\label{eq2}
\theta(u^{+})+\theta(x_{1})=\theta(u^{-})+\theta(x_{0})
\end{equation} since $\theta$ is continuous in $[0,u)$.
Again from the continuity of $\theta$ in $[0,u)$, we have that for each $x,y\in[\theta(u^{+})+\theta(x_{0}),+\infty]$ with $x<y$ there exist two elements $x^{*},y^{*}\in[0,x_{0}]$ satisfying $x=\theta(u^{+})+\theta(x^{*})$ and $y=\theta(u^{+})+\theta(y^{*})$. Hence $\theta(x^{*})<\theta(y^{*})$. This means $x^{*}>y^{*}$ since $\theta$ is non-increasing in $[0,x_{0}]$. Thus $\vartheta(\theta(u^{+})+\theta(x^{*}))\ge\vartheta(\theta(u^{+})+\theta(y^{*}))$ since $O_{\theta,\vartheta}$ is non-decreasing i.e., $\vartheta(x)\ge\vartheta(y)$. Therefore, by the arbitrariness of $x,y$ the last inequality implies that $\vartheta$ is non-increasing in $[\theta(u^{+})+\theta(x_{0}),+\infty]$. On the other hand, from (\ref{eq1}) and (\ref{eq2}), we get that
\begin{equation*}
\vartheta(\theta(x_{0})+\theta(u^{-}))=\vartheta(\theta(x_{0})+\theta(u^{+}))=\vartheta(\theta(x_{1})+\theta(u^{+})),
\end{equation*}
showing that $\vartheta(v)=c$ for any $v\in[\theta(x_{0})+\theta(u^{+}),\theta(x_{1})+\theta(u^{+})]$.

Similarly, there exists an $x_{2}\in[0,x_{0})$ with $x_{2}<x_{1}$ such that
\begin{equation*}
\vartheta(\theta(x_{1})+\theta(u^{-}))=\vartheta(\theta(x_{1})+\theta(u^{+}))=\vartheta(\theta(x_{2})+\theta(u^{+})),
\end{equation*}
showing that $\vartheta(v)=c$ for any $v\in[\theta(x_{1})+\theta(u^{+}),\theta(x_{2})+\theta(u^{+})]$.

Repeating the above process, we get that for arbitrary $x,y\in(0,x_{0}]$, $\vartheta(\theta(x)+\theta(u^{+}))=\vartheta(\theta(y)+\theta(u^{+}))$, i.e., $\vartheta(v)=c$ for any  $v\in[\theta(u^{+})+\theta(x_{0}),+\infty)$. Moreover, if $c=0$, then $\vartheta(\theta(x)+\theta(u^{+}))=0$ with $x\ne0$ and $u\ne0$, contrary to Definition (O2). Hence $c>0$, contrary to the continuity of $O_{\theta,\vartheta}$ since $\vartheta(+\infty)=0$. Therefore, $\theta$ is continuous.
	
	(ii) From Proposition \ref{proposition3.1} (2), if $\theta$ is a monotone function then $\theta$ must be non-increasing. Now, supposing that $\theta$ is not a monotone function, there exist $x_1,x_2,x_3$ such that $\theta$ is non-increasing in $[0,x_{1}]$ and $\theta$ is non-decreasing in $[x_{2},x_{3}]$ in which either $x_{1}=x_{2}$ or $x_{1}<x_{2}$ with $\theta(x_{1})=\theta(x_{2})$. Thus there exists an $x_{4}\in[0,x_{1}]$ such that $\theta(x_{4})=\theta(x_{3})>\theta(x_{1})$ since $\theta(0)=+\infty$. This follows that
\begin{equation*}
\vartheta(\theta(x_{4})+\theta(x))\le\vartheta(\theta(x_{1})+\theta(x))\le\vartheta(\theta(x_{3})+\theta(x))=\vartheta(\theta(x_{4})+\theta(x))
\end{equation*}
 for each $x\in[0,1]$ since $O_{\theta,\vartheta}$ is non-decreasing, i.e., \begin{equation}\label{eq0001}\vartheta(\theta(x_{4})+\theta(x))=\vartheta(\theta(x_{1})+\theta(x))\end{equation} for all $x\in[0,1]$.

On the one hand, there is a $u_{1}\in(0,x_{4})$ such that $\theta(u_{1})>\theta(x_{4})$ since $\theta$ is non-increasing in $[0,x_{1}]$ and, from \eqref{eq0001} we have $\vartheta(\theta(u_{1})+\theta(x_{4}))=\vartheta(\theta(u_{1})+\theta(x_{1}))$. On the other hand, from (i) there exists a $u_{2}\in [0, x_4]$ with $u_{2}<u_{1}$ such that $\theta(u_{2})+\theta(x_{1})=\theta(u_{1})+\theta(x_{4})$ since $\theta(x_{4})>\theta(x_{1})$. Then $\vartheta(\theta(u_{1})+\theta(x_{1}))=\vartheta(\theta(u_{1})+\theta(x_{4}))=\vartheta(\theta(u_{2})+\theta(x_{1}))$.  Obviously, $\vartheta$ is non-increasing in $[2\theta(x_{1}),+\infty]$ since $\theta$ is non-increasing in $[0,x_{1}]$. Thus $\vartheta(v)=c$ for any $v\in[\theta(x_{4})+\theta(x_{1}),\theta(u_{2})+\theta(x_{1})]$.

 Analogously, there exists a $u_{3}\in [0, x_4]$ with $u_{3}<u_{2}$ such that $\theta(u_{3})>\theta(u_{2})$. Thus $\vartheta(v)=c$ for any $v\in[\theta(u_{2})+\theta(x_{1}),\theta(u_{3})+\theta(x_{1})]$.

Repeating the above process, we get that for arbitrary $x,y\in(0,x_{4}]$,
\begin{equation*}
\vartheta(\theta(x)+\theta(x_{1}))=\vartheta(\theta(y)+\theta(x_{1})),
\end{equation*}
i.e., $\vartheta(v)=c$ for any  $v\in[\theta(x_{4})+\theta(x_{1}),+\infty)$ with $c>0$, violating the continuity of $O_{\theta,\vartheta}$ since $\vartheta(+\infty)=0$.

(iii) From Proposition \ref{proposition3.1} (2), it is obviously true.

(iv) Since $\theta$ is non-increasing, $\theta(x)\ge\theta(1)$ for each $x\in[0,1)$. Suppose that there exists a $u\in[0,1)$ such that $\theta(1)=\theta(u)$. Then $O_{\theta,\vartheta}(u,u)=\vartheta(\theta(u)+\theta(u))=\vartheta(\theta(1)+\theta(1))=O_{\theta,\vartheta}(1,1)=1$ with $u\ne1$, contrary to Definition \ref{definition2.1} (O3). Therefore, $\theta(x)>\theta(1)$ for each $x\in[0,1)$.

Conversely, let $\vartheta:[0,+\infty]\rightarrow[-\infty,+\infty]$ be a function satisfying the following three conditions:
\begin{itemize}
  \item[(a)] $\vartheta$ is continuous and non-increasing in $[2\theta(1),+\infty]$;
  \item[(b)] $\vartheta(x)=0$ with $x\in[2\theta(1),+\infty]$ if and only if $x=+\infty$;
  \item[(c)] $\vartheta(x)=1$ with $x\in[2\theta(1),+\infty]$ if and only if $x=2\theta(1)$.
\end{itemize}
Next, we prove that the bivariate function $O_{\theta,\vartheta}:[0,1]^{2}\rightarrow[0,1]$ given by
	$$O_{\theta,\vartheta}(x,y)=\vartheta(\theta(x)+\theta(y))$$
	is an overlap function. The commutativity, continuity and monotonicity are obvious.

Note that for each $x,y\in[0,1]$, $\theta(x)+\theta(y)\in[2\theta(1),+\infty]$ since $\theta$ is non-increasing and continuous. It follows from (b) and (iii) that $\vartheta(\theta(x)+\theta(y))=0$ if and only if $\theta(x)+\theta(y)=+\infty$ if and only if $xy=0$. Moreover, it follows from (c) and  (iv) that $\vartheta(\theta(x)+\theta(y))=1$ if and only if $\theta(x)+\theta(y)=2\theta(1)$ if and only if $xy=1$.

Therefore, by Definition \ref{definition2.1} the bivariate function $O_{\theta,\vartheta}$ is an overlap function.
\end{proof}

The following proposition presents the conditions what $\vartheta$ needs to satisfy when $O_{\theta,\vartheta}(x,y)=\vartheta(\theta(x)+\theta(y))$ is an overlap function.
\begin{proposition}\label{proposition3.2}
	Let $\theta:[0,1]\rightarrow[0,+\infty]$ be a function with both $\theta(x^{-})$ and $\theta(x^{+})$ existing for every $x\in [0,1]$ and $\vartheta:[0,+\infty]\rightarrow[-\infty,+\infty]$ be a function. If the bivariate function $O_{\theta,\vartheta}:[0,1]^{2}\rightarrow[0,1]$ given by
	$$O_{\theta,\vartheta}(x,y)=\vartheta(\theta(x)+\theta(y))$$
is an overlap function then the following three statements hold:
\begin{itemize}
\item[(a)] $\vartheta$ is continuous and non-increasing in $[2\theta(1),+\infty]$;
\item[(b)] $\vartheta(x)=0$ with $x\in[2\theta(1),+\infty]$ if and only if $x=+\infty$;
\item[(c)] $\vartheta(x)=1$ with $x\in[2\theta(1),+\infty]$ if and only if $x=2\theta(1)$.
\end{itemize}
\end{proposition}
\begin{proof}
Suppose that the bivariate function $O_{\theta,\vartheta}$ is an overlap function.
Then from Proposition \ref{proposition03.2} (i) and (ii), $\theta(x)+\theta(y)\in[2\theta(1),+\infty]$ for each $x,y\in[0,1]$.

(a) If $\vartheta$ is not continuous in $[2\theta(1),+\infty]$, then $O_{\theta,\vartheta}(x,y)=\vartheta(\theta(x)+\theta(y))$ is not continuous, a contradiction.

From Proposition \ref{proposition03.2} (ii), $\theta(y)\ge\theta(z)$ for each $x,y,z\in[0,1]$ with $y<z$. Then $\vartheta(\theta(x)+\theta(y))\le\vartheta(\theta(x)+\theta(z))$ since $O_{\theta,\vartheta}$ is non-decreasing. Therefore, $\vartheta$ is non-increasing in $[2\theta(1),+\infty]$.

(b) If $\vartheta(x)=0$ with $x\in[2\theta(1),+\infty]$, then from Proposition \ref{proposition03.2} (i) there exist $y,z\in[0,1]$ such that $x=\theta(y)+\theta(z)$. Hence $\vartheta(\theta(y)+\theta(z))=\vartheta(x)=0$. From Definition \ref{definition2.1} (O2), $yz=0$. Thus from Proposition \ref{proposition03.2} (iii), $x=\theta(y)+\theta(z)=+\infty$.

Conversely, $\vartheta(+\infty)=\vartheta(2\theta(0))=0$ from Proposition \ref{proposition03.2} (iii).

(c) If $\vartheta(x)=1$ with $x\in[2\theta(1),+\infty]$, then from Proposition \ref{proposition03.2} (i) there exist $y,z\in[0,1]$ such that $x=\theta(y)+\theta(z)$. Hence $\vartheta(x)=\vartheta(\theta(y)+\theta(z))=1$. From Definition \ref{definition2.1} (O3), $yz=1$. Thus $x=\theta(y)+\theta(z)=2\theta(1)$.

The converse implication, $\vartheta(x)=\vartheta(2\theta(1))=1$ by Definition \ref{definition2.1} (O3).
\end{proof}

	Notice that from Proposition \ref{proposition3.2}, the values of the function $\vartheta:[0,+\infty]\rightarrow[-\infty,+\infty]$ in $[0,2\theta(1))$ have no influence on the overlap function $O_{\theta,\vartheta}$. Therefore, we only need to consider the values of the function $\vartheta:[0,+\infty]\rightarrow[-\infty,+\infty]$ in $[2\theta(1),+\infty]$ when we investigate the additive generator pair $(\theta,\vartheta)$ of the overlap function $O_{\theta,\vartheta}$.

The following theorem connects Propositions \ref{proposition03.2} with \ref{proposition3.2}.
\begin{theorem}\label{theorem3.2}
	Let $\theta:[0,1]\rightarrow[0,+\infty]$ be a function with both $\theta(x^{-})$ and $\theta(x^{+})$ existing for every $x\in [0,1]$ and $\vartheta:[0,+\infty]\rightarrow[-\infty,+\infty]$ be a function. Then the bivariate function $O_{\theta,\vartheta}:[0,1]^{2}\rightarrow[0,1]$ given by $$O_{\theta,\vartheta}(x,y)=\vartheta(\theta(x)+\theta(y))$$
	is an overlap function if and only if the following statements hold:
	\begin{itemize}
		\item[(1)] $\theta$ is continuous and non-increasing;
		\item[(2)] $\vartheta$ is continuous and non-increasing in $[2\theta(1),+\infty]$;
		\item[(3)] $\theta(x)=+\infty$ if and only if $x=0$;
		\item[(4)] $\vartheta(x)=0$ with $x\in[2\theta(1),+\infty]$ if and only if $x=+\infty$;
	
		\item[(5)]For a given $a\in[0,+\infty)$, the following two conditions are equivalent:
		
		 (i) $\theta(x)=\frac{a}{2}$ if and only if $x=1$;
		
		 (ii) $\vartheta(x)=1$ with $x\in[a,+\infty]$ if and only if $x=a$.
  \end{itemize}	
\end{theorem}
\begin{proof}
Necessity. If the bivariate function $O_{\theta,\vartheta}:[0,1]^{2}\rightarrow[0,1]$ given by $O_{\theta,\vartheta}(x,y)=\vartheta(\theta(x)+\theta(y))$ is an overlap function, then it follows from Propositions \ref{proposition03.2} and \ref{proposition3.2} that the statements (1),(2),(3),(4) and (5) hold.\\

 Sufficiency. In order to prove that the bivariate function $O_{\theta,\vartheta}$ is an overlap function we only need to verify the conditions of Definition \ref{definition2.1}.

 (a) The commutativity, continuity and monotonicity are obvious.

 (b)  According to (3) and (4), $O_{\theta,\vartheta}(x,y)=\vartheta(\theta(x)+\theta(y))=0$ if and only if $\theta(x)+\theta(y)=+\infty$ if and only if $xy=0$.

 (c) From (1) and (5), $O_{\theta,\vartheta}(x,y)=\vartheta(\theta(x)+\theta(y))=1$ if and only if $\theta(x)+\theta(y)=a$ if and only if $xy=1$.

Therefore, $O_{\theta,\vartheta}$ is an overlap function.
\end{proof}

\section{Some transformations of additive generator pairs}\label{section4}
Let $\theta:[0,1]\rightarrow[0,+\infty]$ be a function with both $\theta(x^{-})$ and $\theta(x^{+})$ existing for every $x\in [0,1]$ and $\vartheta:[0,+\infty]\rightarrow[-\infty,+\infty]$ be a function. In this section we first explore some transformations of an additive generator pair $(\theta,\vartheta)$ for the bivariate function $O_{\theta,\vartheta}:[0,1]^{2}\rightarrow[0,1]$ given by
	$$O_{\theta,\vartheta}(x,y)=\vartheta(\theta(x)+\theta(y))$$
	being an overlap function. Then we further describe the pair $(\theta,\vartheta)$ when we extend the functions $\theta$ from $[0,1]$ to $[-\infty,+\infty]$ and $\vartheta$ from $[-\infty,+\infty]$ to $[-\infty,+\infty]$, respectively.

Recall that a function $f:\mathcal{I}\to\overline{\mathbb{R}}$ is called an affine function on the extended real number field, where $\overline{\mathbb{R}}=\mathbb{R} \cup\{+\infty,-\infty\}$ and $\mathcal{I}\subseteq\overline{\mathbb{R}}$, if it can be expressed as
$$f(x) = kx + b$$
with $k,b\in\mathbb{R}$ and by convention,
$$k\cdot(+\infty)+b=+\infty \mbox{ if }k>0\mbox{ and }k\cdot(+\infty)+b=-\infty\mbox{ if }k<0;$$
$$k\cdot(-\infty)+b=-\infty \mbox{ if }k>0\mbox{ and }k\cdot(-\infty)+b=+\infty\mbox{ if }k<0 \mbox{ and}$$
$$f(x)=b\mbox{ for all }x\in\mathcal{I}\mbox{ if }k=0.$$
Then we first have the following proposition.
\begin{proposition}\label{proposition4.1}
	Let $\theta:[0,1]\rightarrow[0,+\infty]$ be a function with both $\theta(x^{-})$ and $\theta(x^{+})$ existing for every $x\in [0,1]$ and $\vartheta:[0,+\infty]\rightarrow[-\infty,+\infty]$ be a function. If $(\theta,\vartheta)$ is an additive generator pair of the overlap function $O_{\theta,\vartheta}:[0,1]^{2}\rightarrow[0,1]$ given by $O_{\theta,\vartheta}(x,y)=\vartheta(\theta(x)+\theta(y))$, then $({h}\circ\theta,\vartheta\circ{g})$ is also an additive generator pair of $O_{\theta,\vartheta}$ where $h:[\theta(1),+\infty]\rightarrow[-\infty,+\infty]$ is given by $h(x)=kx+b$ with certain constants $k,b\in\mathbb{R}$ and $k\ne0$ and, $g:\emph{Ran}(h\circ\theta)+\emph{Ran}(h\circ\theta)\rightarrow[0,+\infty]$ is defined by $g(x)=\frac{x-2b}{k}$ with certain constants $k,b\in\mathbb{R}$ and $k\ne0$.
\end{proposition}
\begin{proof}
We only need to verify that $\vartheta\circ{g}({h}\circ\theta(x)+{h}\circ\theta(y))=\vartheta(\theta(x)+\theta(y))$. Indeed,
\begin{align*}
\vartheta\circ{g}({h}\circ\theta(x)+{h}\circ\theta(y))=&\vartheta\circ{g}(k(\theta(x)+\theta(y))+2b) \\
=&\vartheta(\frac{k(\theta(x)+\theta(y))+2b-2b}{k})\\
=&\vartheta(\theta(x)+\theta(y)).
\end{align*}
\end{proof}

In particular, we have the following corollary.
\begin{corollary}\label{corollary4.1}
	Let $\theta:[0,1]\rightarrow[0,+\infty]$ be a function with both $\theta(x^{-})$ and $\theta(x^{+})$ existing for every $x\in [0,1]$ and $\vartheta:[0,+\infty]\rightarrow[-\infty,+\infty]$ be a function. If $(\theta,\vartheta)$ is an additive generator pair of the overlap function $O_{\theta,\vartheta}:[0,1]^{2}\rightarrow[0,1]$ given by $O_{\theta,\vartheta}(x,y)=\vartheta(\theta(x)+\theta(y))$, then $(\theta^{*},\vartheta^{*})$ is also an additive generator pair of $O_{\theta,\vartheta}$ where $\theta^{*}:[0,1]\rightarrow(-\infty,+\infty]$ is given by $\theta^{*}(x)=\theta(x)+b$ with a certain constant $b\in\mathbb{R}$ and, $\vartheta^{*}:[2\theta(1)+2b,+\infty]\rightarrow[0,1]$ is defined by $\vartheta^{*}(x)=\vartheta(x-2b)$ with a certain constant $b\in\mathbb{R}$.
\end{corollary}

Therefore, in what follows, we only need to consider the case $a=0$, i.e., $\theta(1)=0$ and $\vartheta(0)=1$ whenever applying Theorem \ref{theorem3.2}.

\begin{remark}\emph{In general, in Proposition \ref{proposition4.1} if $h$ is not an affine function then there exist no functions $g$ such that $({h}\circ\theta, \vartheta\circ{g})$ is an  additive generator pair of the overlap function $O_{\theta,\vartheta}$}.\end{remark}
\begin{example}
	\emph{Consider the functions $\theta:[0,1]\rightarrow[0,+\infty]$ and $\vartheta:[0,+\infty]\rightarrow[0,1]$ defined by
		$$\theta(x)=\begin{cases}
			-\mbox{ln}x & \mbox{if}\ x\neq0,\\
			+\infty & \mbox{if}\ x=0
		\end{cases}$$
		and
		$$\vartheta(x)=\begin{cases}
			0        &\mbox{if}\ x=+\infty,\\
			e^{-x} &\mbox{if}\ x\in[0,+\infty),
		\end{cases}$$
		respectively. One can check that $(\theta,\vartheta)$ is an additive generator pair of the overlap function $O_{\theta\,\vartheta}(x,y)=xy$.}

\emph{Now, we consider the function $h:[0,+\infty]\rightarrow[0,+\infty]$ with $h(x)=x^{2}$. If there exists a function $g$ such that $(h\circ\theta,\vartheta\circ{g})$ is an additive generator of $O_{\theta\,\vartheta}=xy$, then $\vartheta\circ{g}(\theta^{2}(x)+\theta^{2}(y))=\vartheta\circ{g}(h\circ\theta(x)+h\circ\theta(y))=\vartheta(\theta(x)+\theta(y))$. Thus for each $x,y\in[0,1]$,
		\begin{equation}\label{equation10}
			g(\theta^{2}(x)+\theta^{2}(y))=\theta(x)+\theta(y)
		\end{equation}
since $\vartheta$ is strictly decreasing and continuous. This follows that $g(({\ln}x)^{2}+({\ln}x)^{2})=-2{\ln}x$, i.e., $g(x)=(2x)^{\frac{1}{2}}$ for any $x\in[0,+\infty]$. Thus $g(\theta^{2}(\frac{1}{2})+\theta^{2}(1))=\sqrt{2}{\ln}2$. On the other hand, $\theta(\frac{1}{2})+\theta(1)={\ln}2$. Therefore, $g(\theta^{2}(\frac{1}{2})+\theta^{2}(1))\neq \theta(\frac{1}{2})+\theta(1)$, contrary to Eq.\eqref{equation10}.}
\end{example}

Below, we cite an important lemma which will be useful in the sequel.
\begin{lemma}[\cite{MK2009}]\label{lemma4.1}
Let \( D \subset \mathbb{R}^N \) be a convex set such that \(\emph{int } D \neq \emptyset\). A function \( f: D \to \mathbb{R} \) is a continuous solution of equation $$f(\frac{x+y}{2})=\frac{f(x)+f(y)}{2}$$ if and only if
\[f(x) = cx + a, \quad x \in D,\]
with certain constants \( c \in \mathbb{R}^N \), \( a \in \mathbb{R} \).
\end{lemma}

\begin{theorem}\label{theorem4.010}
Let $\theta:[0,1]\rightarrow[0,+\infty]$ be a function with both $\theta(x^{-})$ and $\theta(x^{+})$ existing for every $x\in [0,1]$ and $\vartheta:[0,+\infty]\rightarrow[0,1]$ be a strictly decreasing function. Let $(\theta,\vartheta)$ be an additive generator pair of the overlap function $O_{\theta,\vartheta}:[0,1]^{2}\rightarrow[0,1]$ given by $O_{\theta,\vartheta}(x,y)=\vartheta(\theta(x)+\theta(y))$. Then for the respective functions $h:[0,+\infty]\rightarrow[-\infty,+\infty]$ and $g:\emph{Ran}(h\circ\theta)+\emph{Ran}(h\circ\theta)\rightarrow[0,+\infty]$, $({h}\circ\theta,\vartheta\circ{g})$ is also an additive generator pair of $O_{\theta,\vartheta}$ if and only if $h(x)=kx+b$ and $g(x)=\frac{x-2b}{k}$ with certain constants $k,b\in\mathbb{R}$ and $k\ne0$.
\end{theorem}
\begin{proof} Note that from Theorem 3.1, both $\theta$ and $\vartheta$ are continuous and non-increasing because of $(\theta,\vartheta)$ being an additive generator pair of the overlap function $O_{\theta,\vartheta}$, $\theta(1)=0$ and $\vartheta(0)=1$.

Assuming that $({h}\circ\theta,\vartheta\circ{g})$ is an additive generator pair of $O_{\theta,\vartheta}$, we have
\begin{equation*}
	\vartheta(\theta(x)+\theta(y))=\vartheta\circ{g}({h}\circ\theta(x)+{h}\circ\theta(y)).
\end{equation*}
Then $$\theta(x)+\theta(y)={g}({h}\circ\theta(x)+{h}\circ\theta(y))$$ for each $x,y\in[0,1]$ since $\vartheta$ is continuous and strictly decreasing. Letting $u=\theta(x)$ and $v=\theta(y)$, we have
\begin{equation}\label{equation4}
u+v=g(h(u)+h(v))
\end{equation}
for each $u,v\in[0,+\infty]$.

In complete analogy to the proof of Proposition \ref{proposition03.2}, we can prove that {\bf{either}} the following three statements
\begin{itemize}
  \item[(1)] ${h}\circ\theta$ is continuous and non-increasing;
  \item[(2)] ${h}\circ\theta(x)=+\infty$ if and only if $x=0$, i.e., $h(u)=+\infty$ if and only if $u=+\infty$;
  \item[(3)] ${h}\circ\theta(1)\ne-\infty$ and for each $x\in[0,1)$, ${h}\circ\theta(1)<{h}\circ\theta(x)$.
\end{itemize}
{\bf{or}} the following three statements
\begin{itemize}
  \item[(i)] ${h}\circ\theta$ is continuous and non-decreasing;
  \item[(ii)] ${h}\circ\theta(x)=-\infty$ if and only if $x=0$, i.e., $h(u)=-\infty$ if and only if $u=+\infty$;
  \item[(iii)] ${h}\circ\theta(1)\ne+\infty$ and for each $x\in[0,1)$, ${h}\circ\theta(1)>{h}\circ\theta(x)$.
\end{itemize}
hold.

Suppose the first three statements hold. Then ${h}\circ\theta$ is continuous and non-increasing. This means that $h$ is continuous and non-decreasing since $\theta$ is continuous and non-increasing.

Now we prove $h$ is a bijection in $[0,+\infty]$. If $h(a)=h(b)$ with $a,b\in[0,+\infty)$ then from \eqref{equation4} we have $a+v=g(h(a)+h(v))=g(h(b)+h(v))=b+v$ for any $v\in[0,+\infty)$, which implies $a=b$, i.e., $h$ is an injection in $[0,+\infty)$. Thus from (2), $h$ is an injection in $[0,+\infty]$. Consequently, $h$ is a bijection in $[0,+\infty]$ since $h$ is continuous.

Let $t=h(u), s=h(v)$ and $f:[h(0),+\infty]\rightarrow[0,+\infty]$ be defined by $f(x)=h^{-1}(x)$. Then from \eqref{equation4} we get that $$f(t)+f(s)=g(t+s)$$ for each $t,s\in [h(0),+\infty]$. Then $2f(t)=g(2t)$ for each $t\in[h(0),+\infty)$. Hence $\frac{g(2t)}{2}+\frac{g(2s)}{2}=g(\frac{2t+2s}{2})$ for each $t,s\in[h(0),+\infty)$. Since $f$ is continuous in $[h(0),+\infty)$, we have that $g$ is continuous in $[2h(0),+\infty)$. It is easy to see that $[2h(0),+\infty)$ is convex and $\text{int } [2h(0),+\infty) \neq \emptyset$. Then by Lemma \ref{lemma4.1}, $g(x)=cx+m$ for each $x\in[2h(0),+\infty)$ with constants $c,m\in\mathbb{R}$. This follows that $f(t)=\frac{g(2t)}{2}=ct+\frac{m}{2}$ for each $t\in[h(0),+\infty)$. If $c=0$, then $f(x)=\frac{m}{2}$, contrary to the fact that $f$ is a bijection. If $c<0$ then $f$ is strictly decreasing, a contradiction to the fact that $h$ is strictly increasing. Thus $c>0$. On the other hand, if $t=+\infty$ then $f(+\infty)=+\infty=c(+\infty)+\frac{m}{2}$ with $c>0$. Therefore, we finally have $f(t)=ct+\frac{m}{2}$ for each $t\in[h(0),+\infty]$ with a constant $c>0$.

Thus $h(x)=\frac{x-\frac{m}{2}}{c}$ since $f(t)=h^{-1}(t)$ for any $t\in [h(0),+\infty]$. Let $k=\frac{1}{c}$ and $b=-\frac{m}{2c}$. Then $h(x)=kx+b$ for any $x\in[0,+\infty]$ with $k>0$ and $g(x)=\frac{x-2b}{k}$ for any $x\in\mbox{Ran}(h\circ\theta)+\mbox{Ran}(h\circ\theta)$ with $k>0$.

The other being analogous, $h(x)=kx+b$ for any $x\in[0,+\infty]$ with $k<0$ and $g(x)=\frac{x-2b}{k}$ for any $x\in\mbox{Ran}(h\circ\theta)+\mbox{Ran}(h\circ\theta)$ with $k<0$.

Conversely, we only need to verify that $\vartheta\circ{g}({h}\circ\theta(x)+{h}\circ\theta(y))=\vartheta(\theta(x)+\theta(y))$. Indeed,
\begin{align*}
\vartheta\circ{g}({h}\circ\theta(x)+{h}\circ\theta(y))=&\vartheta\circ{g}(k(\theta(x)+\theta(y))+2b) \\
=&\vartheta(\frac{k(\theta(x)+\theta(y))+2b-2b}{k})\\
=&\vartheta(\theta(x)+\theta(y)).
\end{align*}
\end{proof}

Let $h:[0,+\infty]\rightarrow\overline{\mathbb{R}}$ be a function. For any $b\in\mathbb{R}$, define the function $h+b:[0,+\infty]\rightarrow\overline{\mathbb{R}}$ by $$(h+b)(x)=h(x)+b$$ for any $x\in[0,+\infty]$.

\begin{theorem}\label{theorem4.020}
Let $\theta:[0,1]\rightarrow[0,+\infty]$ be a function with both $\theta(x^{-})$ and $\theta(x^{+})$ existing for every $x\in [0,1]$ and $\vartheta:[0,+\infty]\rightarrow[0,1]$ be a strictly decreasing function. Let $(\theta,\vartheta)$ be an additive generator pair of the overlap function $O_{\theta,\vartheta}:[0,1]^{2}\rightarrow[0,1]$ given by $O_{\theta,\vartheta}(x,y)=\vartheta(\theta(x)+\theta(y))$. Then for the respective functions $h:[0,+\infty]\rightarrow[0,+\infty]$ and $g:[0,\vartheta(2h(0))]\rightarrow[0,1]$, $({h}\circ\theta,g\circ\vartheta)$ is also an additive generator pair of $O_{\theta,\vartheta}$ if and only if $h(x)=kx+b$ with certain constants $k,b\in\mathbb{R}$, $k>0$ and $b\ge0$ and, $g(x)=\vartheta\circ{(h+b)^{-1}}\circ\vartheta^{-1}(x)$ with a certain constant $b\in\mathbb{R}$ and $b\ge0$.
\end{theorem}
\begin{proof}Note that from Theorem 3.1, both $\theta$ and $\vartheta$ are continuous and non-increasing because of $(\theta,\vartheta)$ being an additive generator pair of the overlap function $O_{\theta,\vartheta}$, $\theta(1)=0$ and $\vartheta(0)=1$.

Assuming that $({h}\circ\theta,g\circ\vartheta)$ is an additive generator pair of $O_{\theta,\vartheta}$, we have
\begin{equation*}
	\vartheta(\theta(x)+\theta(y))=g\circ\vartheta({h}\circ\theta(x)+{h}\circ\theta(y)).
\end{equation*}
Then $$\theta(x)+\theta(y)=\vartheta^{-1}\circ{g}\circ\vartheta({h}\circ\theta(x)+{h}\circ\theta(y))$$
 for each $x,y\in[0,1]$ since $\vartheta$ is continuous strictly decreasing. Let $u=\theta(x)$, $v=\theta(y)$ and $\overline{g}=\vartheta^{-1}\circ{g}\circ\vartheta$. Thus
\begin{equation}\label{equation5}
u+v=\overline{g}(h(u)+h(v))
\end{equation}
for each $u,v\in[0,+\infty]$.

In complete analogy to the proof of Proposition \ref{proposition03.2}, we can prove that the following three statements hold:
\begin{itemize}
  \item[(i)] ${h}\circ\theta$ is continuous and non-increasing;
  \item[(ii)] ${h}\circ\theta(x)=+\infty$ if and only if $x=0$ i.e. $h(u)=+\infty$ if and only if $u=+\infty$;
  \item[(iii)] for each $x\in[0,1)$, $0\le{h}\circ\theta(1)=h(0)<{h}\circ\theta(x)$.
\end{itemize}

In complete analogy to the proof of Theorem \ref{theorem4.010} we can prove that $h(x)=kx+b$ for any $x\in[0,+\infty]$ with $k>0$ and $b\ge0$ and, $g(x)=\vartheta\circ{(h+b)^{-1}}\circ\vartheta^{-1}(x)$ for each $x\in[0,\vartheta(2h(0))]$ with $k>0$ and $b\ge0$.

Conversely, we only need to verify that $g\circ\vartheta({h}\circ\theta(x)+{h}\circ\theta(y))=\vartheta(\theta(x)+\theta(y))$. Indeed,
\begin{align*}
&{g}\circ\vartheta(\theta\circ{h}(x)+\theta\circ{h}(y))\\
=&{g}\circ\vartheta(k(\theta(x)+\theta(y))+2b)\\
=&\vartheta\circ{(h+b)^{-1}}\circ\vartheta^{-1}\circ\vartheta(k(\theta(x)+\theta(y))+2b)\\
=&\vartheta(\frac{k(\theta(x)+\theta(y))+2b-2b}{k})\\
=&\vartheta(\theta(x)+\theta(y)).
\end{align*}
\end{proof}

\begin{proposition}\label{proposition4.3}
Let $\theta:[0,1]\rightarrow[0,+\infty]$ be a function with both $\theta(x^{-})$ and $\theta(x^{+})$ existing for every $x\in [0,1]$ and $\vartheta:[0,+\infty]\rightarrow[0,1]$ be a function. If $(\theta,\vartheta)$ is an additive generator pair of the overlap function $O_{\theta,\vartheta}:[0,1]^{2}\rightarrow[0,1]$ given by $O_{\theta,\vartheta}(x,y)=\vartheta(\theta(x)+\theta(y))$, then $(\theta\circ{h},\vartheta\circ{g})$ is also an additive generator pair of $O_{\theta,\vartheta}$ when the function $h:[0,1]\rightarrow[0,1]$ satisfies $\theta\circ{h}(x)=c\theta(x)+m$ with certain constants $c,m\in\mathbb{R}$, $c>0$ and $m\ge0$ and, $g:[2m,+\infty]\rightarrow[0,+\infty]$ is defined by $g(x)=\frac{x-2m}{c}$ with certain constants $c,m\in\mathbb{R}$, $c>0$ and $m\ge0$.
\end{proposition}
\begin{proof}
It suffices to check that $\vartheta\circ{g}(\theta\circ{h}(x)+\theta\circ{h}(y))=\vartheta(\theta(x)+\theta(y))$. In fact,
\begin{align*}
	&\vartheta\circ{g}(\theta\circ{h}(x)+\theta\circ{h}(y))\\
	=&\vartheta\circ{g}(c(\theta(x)+\theta(y))+2m)\\
    =&\vartheta(\frac{c(\theta(x)+\theta(y))+2m-2m}{c})\\
	=&\vartheta(\theta(x)+\theta(y)),
\end{align*}
i.e., $(\theta\circ{h},\vartheta\circ{g})$ is an additive generator pair of $O_{\theta,\vartheta}$.
\end{proof}

Generally, the analytical expression of $\theta\circ{h}$ in Proposition \ref{proposition4.3} should be $\theta\circ{h}(x)=c\theta(x)+m$ for any $x\in [0,1]$.
\begin{example}
	\emph{Consider the functions $\theta:[0,1]\rightarrow[0,+\infty]$ and $\vartheta:[0,+\infty]\rightarrow[0,1]$ defined by
		$$\theta(x)=\begin{cases}
			-\mbox{ln}x & \mbox{if}\ x\neq0,\\
			+\infty & \mbox{if}\ x=0
		\end{cases}$$
		and
		$$\vartheta(x)=\begin{cases}
			0        &\mbox{if}\ x=+\infty,\\
			e^{-x} &\mbox{if}\ x\in[0,+\infty),
		\end{cases}$$
		respectively. It is easy to check that $(\theta,\vartheta)$ is an additive generator pair of the overlap function $O_{\theta\,\vartheta}(x,y)=xy$.}

\emph{Now, let the function $h:[0,1]\rightarrow[0,1]$ be defined by $h(x)=\frac{x^{2}+x}{2}$. It is obvious that $\theta\circ{h}(x)=-{\ln}\frac{x^{2}+x}{2}$. If there exists a function $g$ such that $(\theta\circ{h},\vartheta\circ{g})$ is an additive generator of $O_{\theta\,\vartheta}(x,y)=xy$, then $\vartheta\circ{g}(\theta\circ{h}(x)+\theta\circ{h}(y))=\vartheta\circ{g}(-{\ln}\frac{x^{2}+x}{2}-{\ln}\frac{y^{2}+y}{2})=\vartheta(\theta(x)+\theta(y))$. Thus a simple calculation leads to
		\begin{equation}\label{equation11}
			g(-{\ln}\frac{x^{2}+x}{2}-{\ln}\frac{y^{2}+y}{2})=-{\ln}xy
		\end{equation}
for each $x,y\in[0,1]$. Let $x=\frac{1}{2}$ and $y=\frac{1}{2}$. Then we have $g(-{\ln}\frac{3}{4})={\ln}4$. Let $x=\frac{\sqrt{7}-1}{2}$ and $y=1$. Then we have $g(-{\ln}\frac{3}{4})=-{\ln}\frac{\sqrt{7}-1}{2}$, a contradiction.}
\end{example}
\begin{proposition}\label{proposition4.04}
Let $\theta:[0,1]\rightarrow[0,+\infty]$ be a function with both $\theta(x^{-})$ and $\theta(x^{+})$ existing for every $x\in [0,1]$ and $\vartheta:[0,+\infty]\rightarrow[0,1]$ is a strictly decreasing function. If $(\theta,\vartheta)$ is an additive generator pair of the overlap function $O_{\theta,\vartheta}:[0,1]^{2}\rightarrow[0,1]$ given by $O_{\theta,\vartheta}(x,y)=\vartheta(\theta(x)+\theta(y))$, then $(\theta\circ{h},g\circ\vartheta)$ is also an additive generator pair of $O_{\theta,\vartheta}$ when the function $h:[0,1]\rightarrow[0,1]$ satisfies $\theta\circ{h}(x)=c\theta(x)+m$ with certain constants $c,m\in\mathbb{R}$, $c>0$ and $m\ge0$ and, $g:[0,\vartheta(2m)]\rightarrow[0,1]$ is defined by $g(x)=\vartheta\circ{f}^{-1}\circ\vartheta^{-1}(x)$ in which $f:[0,+\infty]\rightarrow[0,+\infty]$ is defined by $f(x)=cx+2m$ with certain constants $c,m\in\mathbb{R}$, $c>0$ and $m\ge0$.
\end{proposition}
\begin{proof}
	 It is enough to verify that $g\circ\vartheta(\theta\circ{h}(x)+\theta\circ{h}(y))=\vartheta(\theta(x)+\theta(y))$. Indeed,
	\begin{align*}
		&{g}\circ\vartheta(\theta\circ{h}(x)+\theta\circ{h}(y))\\
		=&{g}\circ\vartheta(c(\theta(x)+\theta(y))+2m)\\
		=&\vartheta\circ{f^{-1}}\circ\vartheta^{-1}\circ\vartheta(c(\theta(x)+\theta(y))+2m)\\
        =&\vartheta(\frac{c(\theta(x)+\theta(y))+2m-2m}{c})\\
		=&\vartheta(\theta(x)+\theta(y)).
	\end{align*}
	Therefore, $(\theta\circ{h},g\circ\vartheta)$ is an additive generator pair of $O_{\theta,\vartheta}$.
\end{proof}

Notice that we can give an analogous example to show that the analytical expression of $\theta\circ{h}$ in Proposition \ref{proposition4.04} should be $\theta\circ{h}(x)=c\theta(x)+m$ for any $x\in [0,1]$ generally.

From Proposition \ref{proposition4.1} and Theorem \ref{theorem3.2} we can deduce the following theorem.
\begin{theorem}\label{theorem4.1}
	Let $\theta:[0,1]\rightarrow[-\infty,+\infty]$ be a function with both $\theta(x^{-})$ and $\theta(x^{+})$ existing for every $x\in [0,1]$ and $\vartheta:[-\infty,+\infty]\rightarrow[-\infty,+\infty]$ be a function. Then the bivariate function $O_{\theta,\vartheta}:[0,1]^{2}\rightarrow[0,1]$ given by $O_{\theta,\vartheta}(x,y)=\vartheta(\theta(x)+\theta(y))$
	is an overlap function if and only if {\bf{either}} the following five statements
	\begin{itemize}
		\item[(1)] $\theta$ is continuous and non-increasing;
		\item[(2)] $\vartheta$ is continuous and non-increasing in $[2\theta(1),+\infty]$;
		\item[(3)] $\theta(x)=+\infty$ if and only if $x=0$;
		\item[(4)] For $x\in[2\theta(1),+\infty]$, $\vartheta(x)=0$ if and only if $x=+\infty$;
		
		\item[(5)]For a given $a\in(-\infty,+\infty)$, the following two statements are  equivalent:
		
		(i) $\theta(x)=\frac{a}{2}$ if and only if $x=1$;
		
		(ii) For $x\in[a,+\infty]$, $\vartheta(x)=1$ if and only if $x=a$.
	\end{itemize}
	{\bf{or}} the following five statements
	\begin{itemize}
		\item[(1)] $\theta$ is continuous and non-decreasing;
		\item[(2)] $\vartheta$ is continuous and non-decreasing in $[-\infty,2\theta(1)]$;
		\item[(3)] $\theta(x)=-\infty$ if and only if $x=0$;
		\item[(4)] For $x\in[-\infty,2\theta(1)]$, $\vartheta(x)=0$ if and only if $x=-\infty$;
		
		\item[(5)]For a given $a\in(-\infty,+\infty)$, the following two statements are  equivalent:
		
		(i) $\theta(x)=\frac{a}{2}$ if and only if $x=1$;
		
		(ii) For $x\in[-\infty,a]$, $\vartheta(x)=1$ if and only if $x=a$.
	\end{itemize}
hold.
\end{theorem}
\section{Conclusions}\label{section5}
The article mainly supplied the characterization of additive generator pairs of overlap functions (Theorems \ref{theorem4.010}, \ref{theorem4.020} and \ref{theorem4.1}). In particular, the purpose of Theorems \ref{theorem4.010} and \ref{theorem4.020} together with Propositions \ref{proposition4.3} and \ref{proposition4.04} was to give us analytical expressions of some transformations for additive generator pairs. It regrets that we do not know whether the converses of both Propositions \ref{proposition4.3} and \ref{proposition4.04} are true.

\section*{Declaration of competing interest}
The authors declare that they have no known competing financial interests or personal relationships that could have appeared to influence the work reported in this article.

\end{document}